\documentclass[twoside,twocolumn]{article}

  \usepackage{latexsym}
  \usepackage{amsmath}   %needed for \begin{align}... \end{align} environment
  \usepackage{amsfonts}
  \usepackage{amssymb}
  \usepackage{graphicx}
  \usepackage{ifpdf}
  \usepackage{amsthm}
  \usepackage{url}
  \usepackage[noadjust]{cite}
  \usepackage[bottom=3.5cm, top=2.5cm, left=2.5cm, right=2.5cm, columnsep=0.8cm]{geometry}
  \newcommand{\alert}{}

  \newcommand{\field}[1]{\mathbb{#1}}
  
  \newcommand{\R}{\field{R}}

  \newcommand{\HQ}{\field{H}}
  \newcommand{\vect}[1]{\ensuremath{\mbox{\textbf{\textit{#1}}}}}
  \newcommand{\svect}[1]{\ensuremath{\mbox{\textbf{\textit{\small #1}}}}}

  \newcommand{\be}{\begin{equation}}
  \newcommand{\ee}{\end{equation}} % Does not work. Don't know why.
  \newcommand{\bv}[1]{\mathbf{#1}}
  \newcommand{\blue}[1]{{#1}}

  \newcommand{\bfr}{\begin{frame}}
  \newcommand{\efr}{\end{frame}}

\title{Angles between subspaces}

\author{Eckhard Hitzer, Department of Applied Physics, University of Fukui, 910-8507 Japan}
\begin{document}

\DeclareGraphicsExtensions{.pdf}
 
\twocolumn[{\csname @twocolumnfalse\endcsname

\maketitle  % full width title
\thispagestyle{empty}
\pagestyle{empty}

\begin{abstract}
\noindent
We first review the definition of the angle between subspaces and how it is computed using matrix algebra. Then we introduce the Grassmann and Clifford algebra description of subspaces. The geometric product of two subspaces yields the full relative angular information in an explicit manner. We explain and interpret the result of the geometric product of subspaces gaining thus full access to the relative orientation information. 
\vspace{0.5em}


%\vspace{0.5em}

%\subparagraph{AMS Subj. Class.:}
\textbf{AMS Subj. Class.:}
15A66.

\vspace*{1.0\baselineskip}

\end{abstract}
}]

%%%%%%%%%%%%%%%%%%%%%%%%%%%%%%%%%%%%%%%%%%%%
%% MAINMATTER
%%%%%%%%%%%%%%%%%%%%%%%%%%%%%%%%%%%%%%%%%%%%

\section{Introduction}

I first came across Clifford's geometric algebra in the early 90ies in papers on gauge field theory of gravity by J.S.R. Chisholm, struck by the seamlessly compact, elegant, and geometrically well interpretable expressions for elementary particle fields subject to Einstein's gravity. Later I became familiar with D. Hestenes' excellent modern formulation of geometric algebra, which explicitly shows how the geometric product of two vectors encodes their complete relative orientation in the scalar inner product part (cosine) and in the bivector outer product part (sine). 

Geometric algebra can be viewed as an algebra of a vector space and all its subspaces, represented by socalled blades. I therefore often wondered if the geometric product of subspace blades also encodes their complete relative orientation, and how this is done? What is the form of the result, how can it be interpreted and put to further use?
I learned more about this problem, when I worked on the conformal representation of points, point pairs, 
\noindent
\fbox{
\parbox[c]{7.5cm}{Permission to make digital or hard copies of all or part of
this work for personal or classroom use is granted without
fee provided that copies are not made or distributed for
profit or commercial advantage and that copies bear this
notice and the full citation on the first page. To copy
otherwise, or republish, to post on servers or to
redistribute to lists, requires prior specific permission
and/or a fee.}}
lines, planes, circles ans spheres of three dimensional Euclidean geometry \cite{HTBY:Carrier} and was able to find one general formula fully expressing the relative orientation of any two of these objects. Yet L. Dorst (Amsterdam) later asked me if this formula could be generalized to any dimension, because by his experience formulas that work dimension independent are \textit{right}. I had no immediate answer, it seemed to complicated to me, having to deal with too many possible cases. 

But when I prepared for December 2009 a presentation on neural computation and Clifford algebra, I came across a 1983 paper by Per Ake Wedin on angles between subspaces of finite dimensional inner product spaces \cite{PAW:ABS}, which taught me the classical approach. In addition it had a very interesting note on solving the problem, essentially using Grassmann algebra with an additional canonically defined inner product. After that the various bits and pieces came together and began to show the whole picture, the picture which I want to explain in this contribution. 

\section{The angle between two lines} 

\begin{figure}
\begin{center}
  \resizebox{0.25\textwidth}{!}{\includegraphics{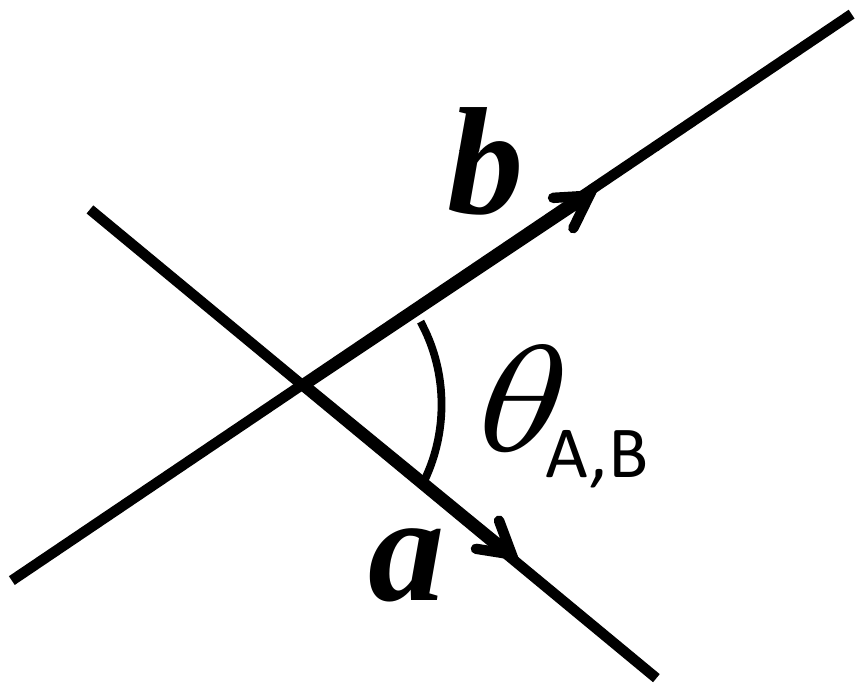}}
  \caption{Angle $\theta_{\mathsf{A},\mathsf{B}}$ between two lines $\mathsf{A}$, $\mathsf{B}$, spanned by unit vectors $\vect{a}$, $\vect{b}$, respectively.
  \label{fg:2lines}}
\end{center}
\end{figure}

To begin with let us look (see Fig. \ref{fg:2lines}) at two lines $\mathsf{A}, \mathsf{B}$ in a vector space $\R^n$, which are spanned by two (unit) vectors $\vect{a},\vect{b} \in \R^n, \vect{a}\cdot\vect{a}= \vect{b}\cdot \vect{b} = 1$: 
\be 
  \mathsf{A} = \text{span}\{\vect{a}\}, \quad
  \mathsf{B} = \text{span}\{\vect{b}\}. 
\ee 
The angle $0\leq \theta_{\mathsf{A},\mathsf{B}}\leq \pi/2$ between lines $\mathsf{A}$ and $\mathsf{B}$ is simply given by
\be 
  \cos \theta_{\mathsf{A},\mathsf{B}} = \vect{a}\cdot\vect{b}.
\ee

\section{Angles between two subspaces (described by principal vectors)}

\begin{figure}
\begin{center}
  \resizebox{0.48\textwidth}{!}{\includegraphics{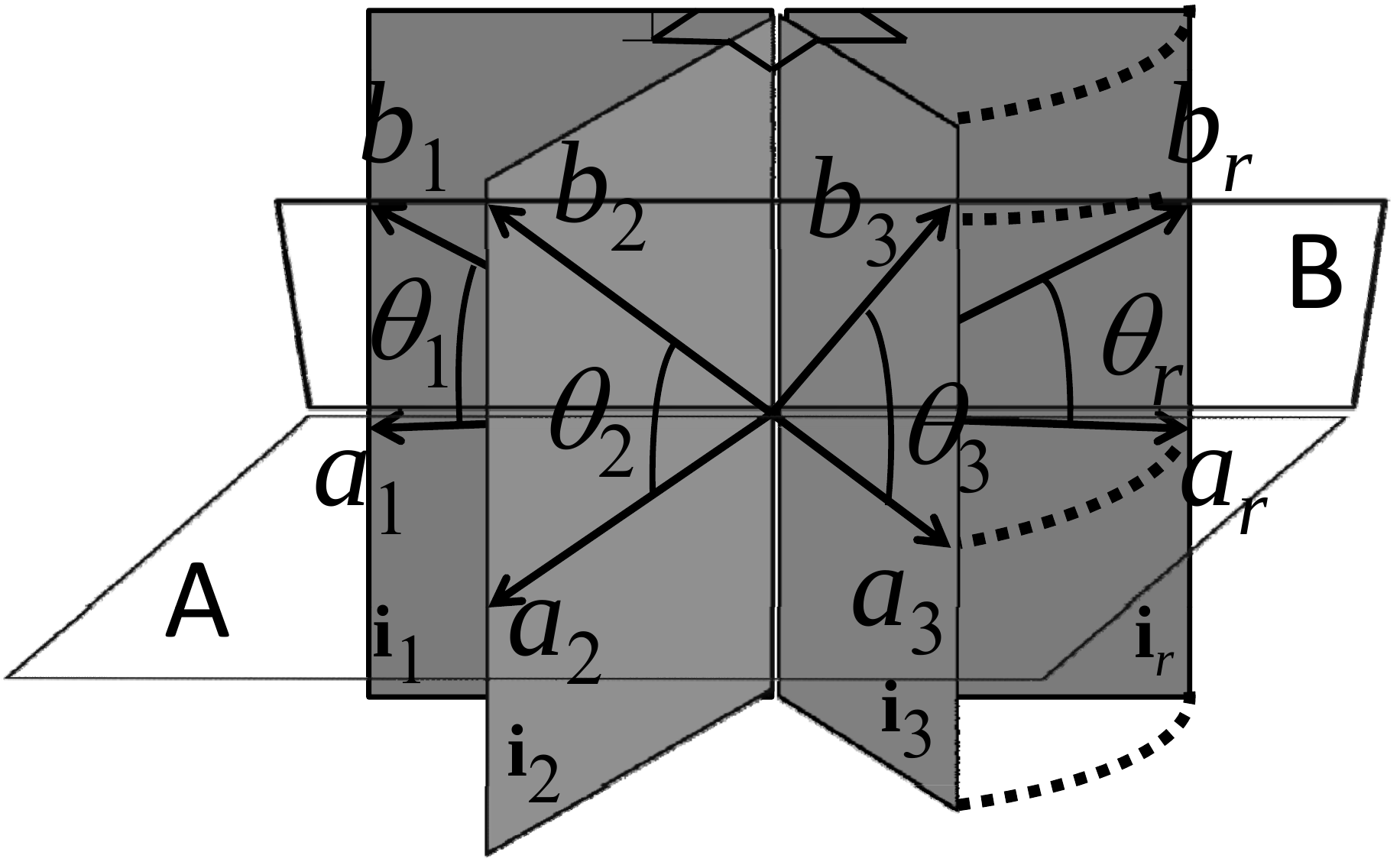}}
  \caption{Angular relationship of two subspaces $\mathsf{A}$, $\mathsf{B}$, spanned by two sets of vectors $\{\vect{a}_1,\ldots,\vect{a}_r\}$, and $\{\vect{b}_1,\ldots,\vect{b}_r\}$, respectively.
  \label{fg:2subspaces}}
\end{center}
\end{figure}

Next let us examine the case of two $r$-dimensional $(r\leq n)$ subspaces $\mathsf{A}, \mathsf{B}$ of an $n$-dimensional Euclidean vector space $\R^n$. The situation is depicted in Fig. \ref{fg:2subspaces}. 
Each subspace $\mathsf{A}$, $\mathsf{B}$ is spanned by a set of $r$ linearly independent vectors
\begin{align} 
  \mathsf{A} &= \text{span}\{\vect{a}_1,\ldots,\vect{a}_r\} \subset \R^n, 
  \nonumber\\
  \mathsf{B} &= \text{span}\{\vect{b}_1,\ldots,\vect{b}_r\} \subset \R^n. 
\end{align}
Using Fig. \ref{fg:2subspaces} we introduce the following notation for principal vectors. The angular relationship between the subspaces $\mathsf{A}$, $\mathsf{B}$ is characterized by a set of $r$ principal angles $\theta_k, 1 \leq k \leq r $, as indicated in Fig. \ref{fg:2subspaces}. A principal angle is the angle between two principal vectors $\vect{a}_k\in \mathsf{A}$ and  $\vect{b}_k\in \mathsf{B}$. The spanning sets of vectors $\{\vect{a}_1,\ldots,\vect{a}_r\}$, and $\{\vect{b}_1,\ldots,\vect{b}_r\}$ can be chosen such that pairs of vectors $\vect{a}_k, \vect{b}_k$ either 
\begin{itemize}
\item agree $\vect{a}_k = \vect{b}_k$, $\theta_k=0$,
\item or enclose a finite angle 
$0 < \theta_{k}\leq \pi/2$. 
\end{itemize}
In addition the pairs of vectors $\{\vect{a}_k, \vect{b}_k\}, 1 \leq k \leq r $ span mutually orthogonal lines (for $\theta_k=0$) and (principal) planes $\bv{i}_k$ (for $0 < \theta_{k}\leq \pi/2$). These mutually orthogonal planes $\bv{i}_k$ are indicated in Fig. \ref{fg:2subspaces}. Therefore if $\vect{a}_k \nparallel \vect{b}_k$ and $\vect{a}_l \nparallel \vect{b}_l$ for $1 \leq k\neq l \leq r$, then plane $\bv{i}_k$ is orthogonal to $\text{plane } \bv{i}_l$. The cosines of the socalled principal angles $\theta_k$ may therefore be $\cos \theta_k = 1$ (for $\vect{a}_k = \vect{b}_k$), or $\cos \theta_k = 0$ (for $\vect{a}_k \perp \vect{b}_k$), or any value $0<\cos \theta_k < 1$. The total angle between the two subspaces $\mathsf{A}$, $\mathsf{B}$ is defined as the product
\be 
  \cos \theta_{\mathsf{A},\mathsf{B}}
  = \cos \theta_1 \cos \theta_2 \ldots \cos \theta_r.
\ee 
In this definition $\cos \theta_{\mathsf{A},\mathsf{B}}$ will automatically be zero if any pair of principal vectors $\{\vect{a}_k, \vect{b}_k\}, 1 \leq k \leq r $ is perpendicular. Then the two subspaces are said to be perpendicular $\mathsf{A} \perp \mathsf{B}$, a familiar notion from three dimensions, where two perpendicular planes $\mathsf{A}$, $\mathsf{B}$ share a common line spanned by $\vect{a}_1 = \vect{b}_1$, and have two mutually orthogonal principal vectors $\vect{a}_2 \perp \vect{b}_2$, which are both in turn orthogonal to the common line vector $\vect{a}_1 $. It is further possible to choose the indexes of the vector pairs $\{\vect{a}_k, \vect{b}_k\}, 1 \leq k \leq r $ such that the principle angles $\theta_k$ appear ordered by magnitude
\be 
  \theta_1 \geq \theta_2 \geq \ldots \geq \theta_r. 
\ee

\section{Matrix algebra computation of angle between subspaces}

The conventional method of computing the angle $\theta_{\mathsf{A},\mathsf{B}}$ between two $r$-dimensional subspaces $\mathsf{A},\mathsf{B} \subset \R^n$ spanned by two sets of vectors $\{\vect{a}'_1,\vect{a}'_2, \ldots \vect{a}'_r \}$ and $\{\vect{b}'_1,\vect{b}'_2, \ldots \vect{b}'_r \}$ is to first arrange these vectors as column vectors into two $n\times r$ matrices
\begin{align} 
  M_\mathsf{A} = [\vect{a}'_1,\ldots,\vect{a}'_r], 
  \quad
  M_\mathsf{B} = [\vect{b}'_1,\ldots,\vect{b}'_r]. 
\end{align}
Then standard matrix algebra methods of QR decomposition and singular value decomposition are applied to obtain 
\begin{itemize}
\item
$r$ pairs of singular unit vectors $\vect{a}_k, \vect{b}_k$ and
\item
$r$ singular values $\sigma_k = \cos \theta_k = \vect{a}_k \cdot \vect{b}_k$. 
\end{itemize}
This approach is very computation intensive.

\section{Even more subtle ways}

Per Ake Wedin in his 1983 contribution \cite{PAW:ABS} to a conference on Matrix Pencils entitled \textit{On Angles between Subspaces of a Finite Dimensional Inner Product Space} first carefully treats the above mentioned matrix algebra approach to computing the angle $\theta_{\mathsf{A},\mathsf{B}}$ in great detail and clarity. Towards the end of his paper he dedicates less than one page to mentioning an alternative method starting out with the words: \textit{But there are even more subtle ways to define angle functions.} 

There he essentially reviews how $r$-dimensional subspaces $\mathsf{A},\mathsf{B} \subset \R^n$ can be represented by $r$-vectors (blades) in Grassmann algebra $A,B \in \Lambda(\R^n)$:
\begin{align} 
  \mathsf{A} = \{\vect{x} \in \R^n | x\wedge A = 0\}, 
  \nonumber \\
  \mathsf{B} =  \{\vect{x} \in \R^n | x\wedge B = 0\}. 
\end{align}
The angle $\theta_{\mathsf{A},\mathsf{B}}$ between the two subspaces $\mathsf{A},\mathsf{B} \in \R^n$ can then be computed in a single step
\be  
  \cos \theta_{\mathsf{A},\mathsf{B}}
  = \frac{A\cdot \widetilde{B}}{|A| |B|}
  = \cos \theta_1 \cos \theta_2 \ldots \cos \theta_r, 
\ee
where the inner product is canonically defined on the Grassmann algebra $\Lambda(\R^n)$ corresponding to the geometry of $\R^n$. The tilde operation is the reverse operation representing a dimension dependent sign change
$\widetilde{B} = (-1)^{\frac{r(r-1)}{2}}B$,
and $|A|$ represents the norm of blade $A$, i.e. 
$|A|^2 = A \cdot \widetilde{A}$, and similarly $|B|^2 = B \cdot \widetilde{B}$. Wedin refers to earlier works of L. Andersson \cite{LA:ABS} in 1980, and a 1963 paper of Q.K. Lu \cite{QKL:EGoES}. 

Yet equipping a Grassmann algebra $\Lambda(\R^n)$ with a canonical inner product comes close to introducing Clifford's geometric algebra $Cl_n = Cl(\R^n)$. And there is another good reason to do that, as e.g. H. Li explains in his excellent 2008 textbook \textit{Invariant Algebras and Geometric Reasoning} \cite{HL:IAaGR}: \textit{... to allow sums of angles to be advanced invariants, the inner-product Grassmann algebra must be extended to the Clifford algebra ...} This is why I have decided to immediately begin in the next section with Clifford's geometric algebra instead of first reviewing inner-product Grassmann algebra.

\section{Clifford (geometric) algebra}

\begin{figure}
\begin{center}
  \resizebox{0.15\textwidth}{!}{\includegraphics{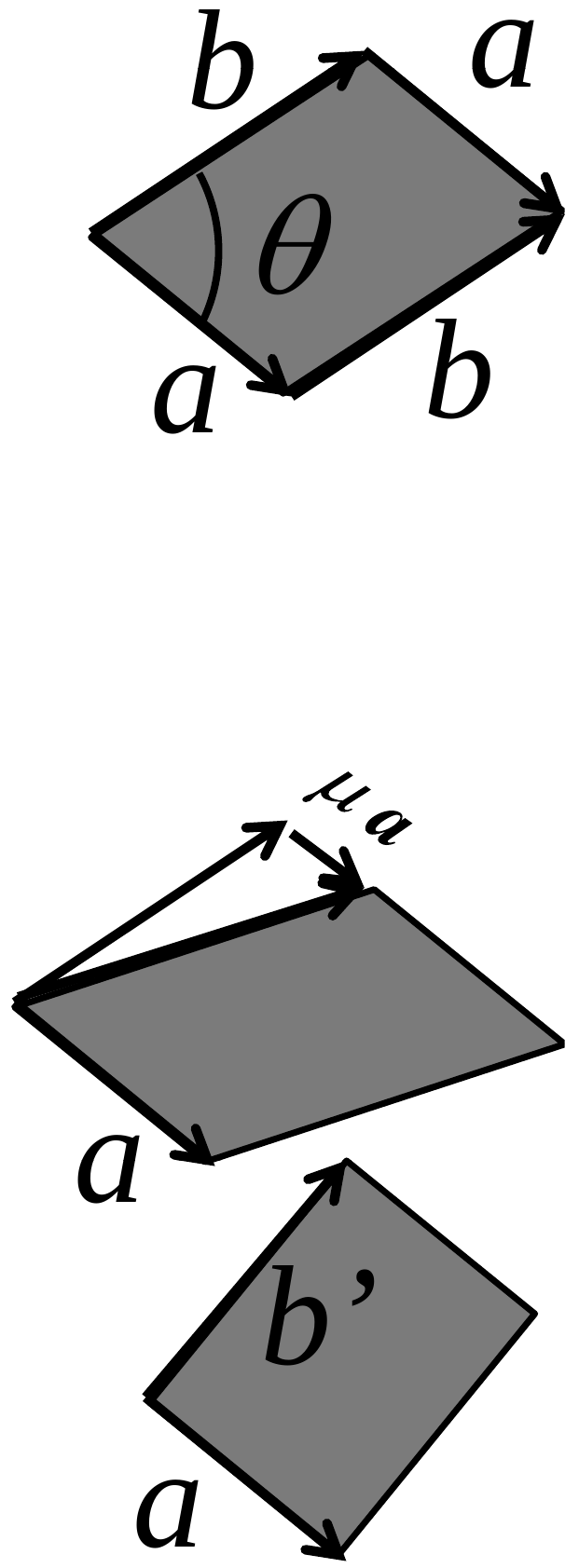}}
  \caption{Bivectors $\vect{a}\wedge \vect{b}$ as oriented area elements can be reshaped (e.g. by $\vect{b} \rightarrow \vect{b}+\mu \vect{a}, \mu \in \R$) without changing their value (area and orientation). The bottom figure shows orthogonal reshaping into the form of an oriented rectangle. 
  \label{fg:bivectors}}
\end{center}
\end{figure}

Clifford (geometric) algebra is based on the geometric product of vectors $\vect{a},\vect{b} \in \R^{p,q}, p+q=n$ 
\begin{equation}
  \vect{a}\vect{b} = \vect{a}\cdot\vect{b} + \vect{a}\wedge\vect{b},
\end{equation}
and the associative algebra $Cl_{p,q}$ thus generated with $\R$ and $\R^{p,q}$ as subspaces of $Cl_{p,q}$. $\vect{a}\cdot\vect{b}$ is the symmetric inner product of vectors and $ \vect{a}\wedge\vect{b}$ is Grassmann's outer product of vectors representing the oriented parallelogram area spanned by $\vect{a},\vect{b}$, compare Fig. \ref{fg:bivectors}.

As an example we take the Clifford geometric algebra $Cl_{3}=Cl_{3,0}$ of three-dimensional (3D) Euclidean space $\R^3=\R^{3,0}$.
$\R^3$ has an orthonormal basis $\{\bv{e}_1, \bv{e}_2, \bv{e}_3\}$. 
$Cl_{3}$ then has an eight-dimensional basis of  
\be
  \label{eq:G3basis}
  \{{1}, 
    \underbrace{\bv{e}_1, \bv{e}_2, \bv{e}_3}_{\text{vectors}},
    {\underbrace{\bv{e}_2\bv{e}_3, \bv{e}_3\bv{e}_1, \bv{e}_1\bv{e}_2}_{\text{area bivectors}}},  
    \underbrace{i=\bv{e}_1\bv{e}_2\bv{e}_3}_{\text{volume trivector}}\}.
\ee
Here $i$ denotes the unit trivector, i.e. the oriented volume of a unit cube, with $i^2=-1$. 
The even grade subalgebra $Cl_{3}^+$ is
isomorphic to Hamilton's quaternions $\HQ$. 
Therefore elements of $Cl_{3}^+$ are also called \textit{rotors} (rotation operators), 
rotating vectors and multivectors of $Cl_{3}$. 

In general $Cl_{p,q}, p+q=n$ is composed of so-called $r$-vector subspaces spanned by the induced bases
\be 
  \label{eq:rvecbasis}
  \{\vect{e}_{k_1} \vect{e}_{k_2} \ldots  \vect{e}_{k_r}
   \mid 1 \leq k_1 < k_2 < \ldots < k_r \leq n \},
\ee
each with dimension $\binom{r}{n}$. The total dimension of the $Cl_{p,q}$ therefore becomes $\sum_{r=0}^n \binom{r}{n} = 2^n$.

General elements called \textit{multivectors}
  $M \in Cl_{p,q}, p+q=n,$ have $k$-vector parts ($0\leq k \leq n$):
  \alert{scalar} part
  $Sc(M) = \langle M \rangle = \langle M \rangle_0 = M_0 \in \R$, 
  \alert{vector} part
  $\langle M \rangle_1 \in \R^{p,q}$, 
  \alert{bi-vector} part
  $\langle M \rangle_2$,  \ldots, 
  and
  \alert{pseudoscalar} part $\langle M \rangle_n\in\bigwedge^n\R^{p,q}$
\begin{equation}\label{eq:MVgrades}
    M  =  \sum_{A=1}^{2^n} M_{A} \vect{e}_{A}
       =  \langle M \rangle + \langle M \rangle_1 + \langle M \rangle_2 + \ldots +\langle M \rangle_n \, .
\end{equation}

The \textit{reverse} of $M \in Cl_{p,q}$ defined as
\begin{equation}\label{eq:MVrev}
  \widetilde{M}=\; \sum_{k=0}^{n}(-1)^{\frac{k(k-1)}{2}}\langle M \rangle_k,
\end{equation}
often replaces \blue{complex conjugation and quaternion conjugation}. Taking the reverse is equivalent to reversing the order of products ob basis vectors in the basis blades of \eqref{eq:rvecbasis}. 
For example the reverse of the bivector $\vect{e}_1\vect{e}_2$ is
\be 
  \widetilde{\vect{e}_1\vect{e}_2}
  = \vect{e}_2\vect{e}_1 = -\vect{e}_1\vect{e}_2,
\ee 
because only the antisymmetric outer product part 
$\vect{e}_2\vect{e}_1 = \vect{e}_2\wedge\vect{e}_1 = - \vect{e}_1\wedge \vect{e}_2$ is relevant.

The \alert{scalar product} of two multivectors $M, \widetilde{N} \in Cl_{p,q}$ is defined as
\be
    M \ast \widetilde{N} 
    = \langle M\widetilde{N} \rangle 
    = \langle M\widetilde{N} \rangle_0.
\ee
  For $M, \widetilde{N} \in Cl_{n}=Cl_{n,0}$ we get $M\ast \widetilde{N}=\sum_{A} M_A N_A.$
  The \blue{modulus} $|M|$ of a multivector $M \in Cl_{n}$ is defined as 
  \be
     |M|^2 = {M\ast\widetilde{M}}= {\sum_{A} M_A^2}.
  \ee

\subsection{Subspaces described in geometric algebra}

In $Cl_n$ symmetric inner product part of two vectors 
$\vect{a}=a_1\vect{e}_1 + a_2\vect{e}_2$, $\vect{b}=b_1\vect{e}_1 + b_2\vect{e}_2$ yields the expected result
\be 
  \vect{a}\cdot \vect{b}
  = a_1b_1 + a_2b_2 
  = |\vect{a}| |\vect{b}| \cos \theta_{\svect{a},\svect{b}}. 
\ee 
Whereas the antisymmetric outer product part gives the bivector, which represents the oriented area of the parallelogram spanned by $\vect{a}$ and $\vect{b}$
\be 
  \vect{a}\wedge \vect{b}
  = (a_1b_2 - a_2b_1) \vect{e}_1\vect{e}_2
  = |\vect{a}| |\vect{b}| \sin \theta_{\svect{a},\svect{b}}\vect{e}_1\vect{e}_2.
\ee 
The parallelogram has the (signed) scalar area $|\vect{a}| |\vect{b}| \sin \theta_{\svect{a},\svect{b}}$ and its orientation in the space $\R^n$ is given by the oriented unit area bivector $\vect{e}_1\vect{e}_2$. 

Two non-zero vectors $\vect{a}$ and $\vect{b}$ are parallel, if and only if $\vect{a}\wedge \vect{b} = 0$, i.e. if and only if $\sin \theta_{\svect{a},\svect{b}} = 0$
\be 
  \vect{a}\wedge \vect{b} = 0
  \Leftrightarrow
  \vect{a} \parallel \vect{b} 
  \Leftrightarrow
  \vect{b} = \alpha \vect{a}, \alpha \in \R.
\ee  
We can therefore use the outer product to represent a line $\mathsf{A}=\text{span}\{\vect{a}\}$ with direction vector $\vect{a}\in\R^n$ as
\be 
  \mathsf{A} = \{\vect{x}\in \R^n \mid \vect{x}\wedge \vect{a}=0\}.
\ee 
Moreover, bivectors can be freely reshaped (see Fig. \ref{fg:bivectors}), e.g.
\be 
  \vect{a}\wedge \vect{b} = \vect{a}\wedge (\vect{b} +\mu \vect{a}), \mu \in \R,
\ee 
because due to the antisymmetry $\vect{a}\wedge \vect{a} = 0$. This reshaping allows to (orthogonally) reshape a bivector to rectangular shape
\be 
  \vect{a}\wedge \vect{b} = \vect{a} \vect{b}', 
  \vect{a} \perp \vect{b} (\text{i.e. } \vect{a} \cdot \vect{b}'=0)
\ee 
as indicated in Fig. \ref{fg:bivectors}. The shape may even chosen as square or circular, depending on the application in mind. 

The total antisymmetry of the trivector 
$\vect{x}\wedge\vect{a}\wedge \vect{b}$ 
means that 
\be 
  \vect{x}\wedge\vect{a}\wedge \vect{b} = 0
  \Leftrightarrow
  \vect{x} = \alpha \vect{a}+\beta\vect{b}, 
  \alpha,\beta  \in \R.
\ee
Therefore a plane $\mathsf{B}$ is given by a simple bivector (also called 2-blade) $B = \vect{a}\wedge \vect{b}$ as
\be 
   \mathsf{B} = \{\vect{x}\in \R^n \mid \vect{x}\wedge B=0\}.
\ee 

A three-dimensional volume subspace $\mathsf{C}$ is similarly given by a 3-blade
$C = \vect{a}\wedge \vect{b}\wedge \vect{c}$ as
\be 
   \mathsf{C} = \{\vect{x}\in \R^n \mid \vect{x}\wedge C=0\}.
\ee 

Finally a \blue{blade} $D_r=\vect{b}_1\wedge\vect{b}_2\wedge\ldots\wedge\vect{b}_r, \vect{b}_l\in\R^{n}, 1\leq l \leq r \leq n$ describes an $r$-dimensional vector \textit{subspace} 
  \be
    \mathsf{D}=\{ \vect{x}\in \R^{p,q} | \vect{x}\wedge D =0 \}.
  \ee 
  Its \blue{dual blade} 
  \be
  D^{\ast}= Di_n^{-1}
  \ee
  describes the \textit{complimentary} $(n-r)$-dimensional vector subspace $\mathsf{D}^{\perp}$. The magnitude of the blade $D_r \in Cl_n$ is nothing but the volume of the $r$-dimensional parallelepiped spanned by the vectors 
$\{\vect{b}_1,\vect{b}_2,\ldots,\vect{b}_r\}$. 

Just as we were able to orthogonally reshaped a bivector to rectangular or square shape we can reshape every $r$-blade $A_r$ to a geometric product of mutually orthogonal vectors 
\be 
  A_r = \vect{a}'_1\wedge \vect{a}'_2\wedge \ldots \vect{a}'_r
  = \vect{a}_1 \vect{a}_2 \ldots \vect{a}_r,
\ee 
with pairwise orthogonal and anticommuting vectors $\vect{a}_1 \perp \vect{a}_2 \perp \ldots \perp \vect{a}_r$. The reverse $\widetilde{A}_r$ of the geometric product of orthogonal vectors $A_r=\vect{a}_1 \vect{a}_2 \ldots \vect{a}_r$ is therefore clearly
\be 
  (\vect{a}_1 \vect{a}_2 \ldots \vect{a}_r)^{\sim}
  = \vect{a}_r \ldots  \vect{a}_2 \vect{a}_1
  = (-1)^{\frac{r(r-1)}{2}} \vect{a}_1 \vect{a}_2 \ldots \vect{a}_r,
\ee 
by simply counting the number $\frac{r(r-1)}{2}$ of permutations necessary.

Paying attention to the dimensions we find that the outer product of an $r$-blade $B$ with a vector $\vect{a}$ increases the dimension (grade) by $+1$
\be 
  \vect{a}\wedge B = \langle\vect{a} B \rangle_{r+1}. 
\ee 
Opposite to that, the inner product (or left contraction) with a vector lowers the dimension (grade) by $-1$
\be 
  \vect{a}\cdot B = \langle\vect{a} B \rangle_{r-1}. 
\ee
The geometric product of two $r$-blades $A,B$ contains therefore at most the following grades
\be 
  AB =
  \langle AB\rangle_0 + \langle AB\rangle_2 + \ldots 
  +\langle AB\rangle_{2\text{min}(r,[n/2])},
\ee 
where the limit $[n/2]$ (entire part of $n/2$) is due to the dimension limit of $\R^n$. 

The inner product of vectors is properly generalized in geometric algebra by introducing the (left) contraction of the $r$-blade $A=A_r$ onto the $s$-blade $B=B_s$ as
\be 
  A_r\rfloor B_s = \langle A B \rangle_{s-r}. 
\ee 
For blades of equal grade ($r=s$) we thus get the symmetric scalar
\be 
  A_r\rfloor B_r = \langle A B \rangle_0 = \langle BA \rangle_0 = A \ast B. 
\ee 

Finally the product of a blade with its own reverse is necessarily scalar. Introducing orthogonal reshaping this scalar is seen to be
\be 
  A_r\widetilde{A}_r
  = \vect{a}_1 \ldots \vect{a}_r\vect{a}_r\ldots\vect{a}_1
  = \vect{a}_1^2 \ldots \vect{a}_r^2 = |A_r|^2,
\ee 
therefore
\be 
  |A_r| = |\vect{a}_1|\ldots|\vect{a}_r|.
\ee 

Every $r$-blade $A_r$ can therefore be written as a product of the scalar magnitude $|A_r| $ times the geometric product of exactly $r$ mutually orthogonal unit vectors $\{\widehat{\vect{a}}_1, \ldots, \widehat{\vect{a}}_r\}$
\be 
  A_r = |A_r| \widehat{\vect{a}}_1 \widehat{\vect{a}}_2 \ldots \widehat{\vect{a}}_r.
\ee 
Please note well, that this \textit{rewriting} of an $r$-blade in geometric algebra does not influence the overall result on the left side, the $r$-blade $A_r$ is before and after the rewriting the very same element of the geometric algebra $Cl_n$. But for the geometric interpretation of the geometric product $AB$ of two $r$-blades $A,B\in Cl_n$ the orthogonal reshaping is indeed a key step. 

After a short discussion of reflections and rotations implemented in geometric algebra, we return to the geometric product of two $r$-blades $A,B\in Cl_n$ and present our key insight.

\subsection{Reflections and rotations}

\begin{figure}
\begin{center}
  \resizebox{0.2\textwidth}{!}{\includegraphics{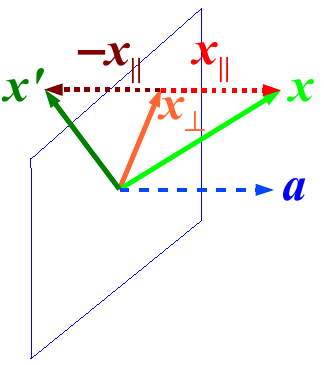}}
  \resizebox{0.2\textwidth}{!}{\includegraphics{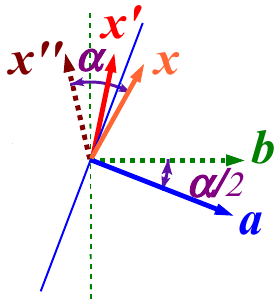}}
  \caption{Reflection at a plane with normal $\vect{a}$ (left) and rotation as double reflection at planes normal to $\vect{a}$, $\vect{b}$ (right). 
  \label{fg:refrot}}
\end{center}
\end{figure}

A simple application of the geometric product is shown in Fig. \ref{fg:refrot} (left) to the reflection of a point vector $\vect{x}$ at a plane with normal vector $\vect{a}$, which means to reverse the component of $\vect{x}$ parallel to $\vect{a}$ (perpendicular to the plane)
\be 
  \vect{x} \longrightarrow \vect{x}' = -\vect{a}^{-1}\vect{x}\vect{a},
  \vect{a}^{-1} = \frac{\vect{a}}{\vect{a}^2}. 
\ee 
Two reflections lead to a rotation by twice the angle between the reflection planes as shown in Fig. \ref{fg:refrot} (right)
\be 
  \vect{x} \longrightarrow \vect{x}'' = \vect{a}^{-1}\vect{b}^{-1}\vect{x}\vect{a}\vect{b}
  = (\vect{a}\vect{b})\vect{x}\vect{a}\vect{b}
  = R^{-1} \vect{x} R, 
\ee 
with rotation operator (rotor) 
$R=\vect{a}\vect{b} \propto \cos \theta_{\svect{a},\svect{b}}+\bv{i}_{\svect{a},\svect{b}}\sin \theta_{\svect{a},\svect{b}}$, where the unit bivector $\bv{i}_{\svect{a},\svect{b}}$ represents the plane of rotation. 

\section{Geometric information in the geometric product of two subspace $r$-blades}

From the foregoing discussion of the representation of $r$-dimensional subspaces $\mathsf{A}, \mathsf{B}\subset \R^n$ by the blades $A=\vect{a}'_1\wedge \ldots \vect{a}'_r$ and $B=\vect{b}'_1\wedge \ldots \vect{b}'_r$, from the freedom of orthogonally reshaping these blades and factoring out the blade magnitudes $|A|$ and $|B|$, and from the classical results of matrix algebra, we now know that we can rewrite the geometric product $AB$ in mutually orthogonal products of pairs of principal vectors $\vect{a}_k,\vect{b}_k, 1\leq k \leq r$
\begin{align} \label{eq:ABrots}
  A\widetilde{B}
  = \vect{a}_1\vect{a}_2\ldots \vect{a}_r\vect{b}_r\ldots \vect{b}_2\vect{b}_1
  = \vect{a}_1\vect{b}_1\vect{a}_2\vect{b}_2\ldots \vect{a}_r\vect{b}_r.
\end{align}
The geometric product 
\begin{align} 
\vect{a}_r\vect{b}_r
  &= |\vect{a}_r||\vect{b}_r|(\cos \theta_{\svect{a}_r,\svect{b}_r} +\bv{i}_{\svect{a}_r,\svect{b}_r}\sin\theta_{\svect{a}_r,\svect{b}_r})
  \nonumber \\
  &= |\vect{a}_r||\vect{b}_r| (c_r +\bv{i}_rs_r), 
\end{align}
with $c_r=\cos \theta_{\svect{a}_r,\svect{b}_r}$ and
 $s_r=\sin \theta_{\svect{a}_r,\svect{b}_r}$
in the above expression for  $A\widetilde{B}$ is composed of a scalar and a bivector part. The latter is proportional to the unit bivector $\bv{i}_r$ representing a (principal) plane orthogonal to \textit{all} other principal vectors. $\bv{i}_r$ therefore commutes with all other principal vectors, and hence the whole product $\vect{a}_r\vect{b}_r$ (a rotor) commutes with all other principal vectors. A completely analogous consideration applies to all products of pairs of principal vectors, which proofs the second equality in \eqref{eq:ABrots}. 

We thus find that we can always rewrite the product $A\widetilde{B}$ as a product of rotors
\begin{align} 
  &A\widetilde{B}
  \nonumber \\
  &= |A||B|(c_1 +\bv{i}_1s_1)(c_2 +\bv{i}_2s_2)\ldots(c_r +\bv{i}_rs_r)
  \nonumber \\
  &= |A||B|(c_1c_2\ldots c_r +
  \nonumber \\
  &\quad+s_1c_2\ldots c_r\bv{i}_1+c_1s_2\ldots c_r\bv{i}_2+\ldots+c_1c_2\ldots s_r\bv{i}_r+
  \nonumber \\
  &\quad \vdots
  \nonumber \\
  &\quad+s_1s_2\ldots s_r\bv{i}_1\bv{i}_2\ldots\bv{i}_r)
\end{align}

We realize how the scalar part 
$\langle A\widetilde{B} \rangle_0=|A||B|c_1c_2\ldots c_r$, 
the bivector part 
$\langle A\widetilde{B} \rangle_2=|A||B|(s_1c_2\ldots c_r\bv{i}_1+c_1s_2\ldots c_r\bv{i}_2+\ldots+c_1c_2\ldots s_r\bv{i}_r)$, 
etc., up to the $2r$-vector (or $2[n/2]$-vector) part 
$\langle A\widetilde{B} \rangle_{2r}=|A||B|s_1s_2\ldots s_r\bv{i}_1\bv{i}_2\ldots\bv{i}_r$ of the geometric product $A\widetilde{B}$ arise and what information they carry. 

Obviously the scalar part yields the cosine of the angle between the subspaces represented by the two $r$-vectors $A,B\in Cl_n$
\be 
  \cos \theta_{\mathsf{A}\mathsf{B}} 
  = \frac{\langle A\widetilde{B} \rangle_0}{|A||B|}
  = \frac{A\ast\widetilde{B}}{|A||B|},
\ee 
which exactly corresponds to P.A. Wedin's formula from inner-product Grassmann algebra. 

The bivector part consists of a sum of (principal) plane bivectors, which can in general be uniquely decomposed  into its constituent
sum of 2-blades by the method of Riesz, described also in \cite{HS:CAtoGC},  chapter 3-4, equation (4.11a) and following. 

The magnitude of the $2r$-vector part allows to compute the product of all sines of principal angles
\be 
  s_1s_2\ldots s_r = \pm \frac{|\langle A\widetilde{B} \rangle_{2r}|}{|A||B|}. 
\ee

Let us finally refine our considerations to two general $r$-dimensional subspaces $\mathsf{A},\mathsf{B}$, which we take to partly intersect and to be partly perpendicular. We mean by that, that the dimension of the intersecting subspace be $s\leq r$ ($s$ is therefore the number of principal angles equal zero), and the number of principle angles with value $\pi/2$ be $t\leq r-s$. For simplicity we work with normed blades (i.e. after dividing with 
$|A||B|$. The geometric product of the the $r$-blades $A,B\in Cl_n$ then takes the form
\begin{align} 
  &A\widetilde{B}
  \nonumber \\
  &= (c_{s+1}c_{s+2}\ldots c_{r-t} +
  \nonumber \\
  &\quad+s_{s+1}c_{s+2}\ldots c_{r-t}\bv{i}_{s+1}
        +c_{s+1}s_{s+2}\ldots c_{r-t}\bv{i}_{s+2}+
  \nonumber \\
  &\quad\ldots
        +c_{s+1}c_{s+2}\ldots s_{r-t}\bv{i}_{r-t}+
  \nonumber \\
  &\quad \vdots
  \nonumber \\
  &\quad+s_{s+1}s_{s+2}\ldots s_{r-t}\bv{i}_{s+1}\bv{i}_{s+2}\ldots\bv{i}_{r-t})
  \bv{i}_{r-t+1} \ldots \bv{i}_{r}.
\end{align}
We thus see, that apart from the integer dimensions $s$ for parallelity (identical to the dimension of the meet of blade $A$ with blade $B$) and $t$ for perpendicularity, the lowest non-zero grade of dimension $2t$ gives the relevant angular measure
\be 
  \cos\theta_{\mathsf{A}\mathsf{B}}=\cos \theta_{s+1}\cos \theta_{s+2}\ldots \cos \theta_{r-t}.
\ee 
While the maximum grade part gives again the product of the corresponding sinuses
\be 
  \sin \theta_{s+1}\sin \theta_{s+2}\ldots \sin \theta_{r-t}.
\ee

Dividing the product $A\widetilde{B}$ by its lowest grade part  
$c_{s+1}c_{s+2}\ldots c_{r-t} \bv{i}_{r-t+1} \ldots \bv{i}_{r}$ gives a multivector with maximum grade $2(r-t-s)$, scalar part one, and bivector part 
\be 
t_{s+1}\bv{i}_{s+1}
 +t_{s+2} \bv{i}_{s+2}+\ldots
 +t_{r-t}\bv{i}_{r-t},
\ee
where $t_k = \tan \theta_k $. 
Splitting this bivector into its constituent bivector parts further yields the (principal) plane bivectors and the tangens values of the principle angles $\theta_k, s<k\leq r-t$. This is the only somewhat time intensive step.

\section{Conclusion}

Let us conclude by discussing possible future applications of these results. The complete relative orientation information in $A\widetilde{B}$ should be ideal for a subspace structure self organizing map (SOM) type of neural network. Not only data points, but the topology of whole data subspace structures can then be faithfully mapped to lower dimensions. Our discussion gives meaningful results for partly intersecting and partly perpendicular subspaces. Apart from extracting the bivector components, all computations are done by multiplication. Projects like fast Clifford algebra hardware developed at the TU Darmstadt (D. Hildenbrand et al) should be of interest for applying the results of the paper to high dimensional data sets. An extension to offset subspaces (of projective geometry) and $r$-spheres (of conformal geometric algebra) may be possible.

\section*{Acknowledgments}
  Soli deo gloria. I do thank my dear family, H. Ishi, D. Hildenbrand and V. Skala.

\bibliographystyle{plain}

\begin{thebibliography}{99}

\bibitem{LA:ABS}
L. Andersson, The concepts of angle between subspaces ... unpublished notes, Umea (1980).

\bibitem{PAW:ABS}
P. A. Wedin, On angles between subspaces of a finite-dimensional inner product space, in Matrix Pencils, Bo Kagstram and Axel Ruhe, eds., Springer-Verlag, Berlin, 1983, pp. 263--285.

\bibitem{LD:GAfCS}
L. Dorst et al, Geometric Algebra for Comp. Sc., Morgan Kaufmann, 2007.

\bibitem{HL:IAaGR}
H. Li, Invariant Algebras and Geometric Reasoning, World Scientific, Singapore, 2008.

\bibitem{QKL:EGoES}
Q.K. Lu, The elliptic geometry of extended spaces, Acta Math. Sinica, 13 (1963), pp. 49--62; translated as Chinese Math. 4 (1963), pp. 54--69.

\bibitem{HS:CAtoGC}
D. Hestenes, G. Sobczyk, Clifford Algebra to Geometric Calculus, Kluwer, 1984.

\bibitem{HTBY:Carrier}
E. Hitzer, K. Tachibana, S. Buchholz, I. Yu
Carrier method for the general evaluation and control of pose, molecular conformation, tracking, and the like
Advances in Applied Clifford Algebras, Vol. 19(2), pp. 339-364 (2009). 


 
    
\end{thebibliography}

\end{document}